\documentclass[a4paper]{amsart}
\usepackage{times, url}
\usepackage[usenames]{color}
\usepackage{mathrsfs}
\usepackage{enumerate}
\usepackage{setspace}

\newcommand{\<}{\begin{equation}}
\newcommand{\?}{\end{equation}}

\newcommand{\Q}{\mathbb{Q}}
\newcommand{\R}{\mathbb{R}}
\newcommand{\Z}{\mathbb{Z}}
\newcommand{\N}{\mathbb{N}}
\newcommand{\C}{\mathbb{C}}

\newcommand{\cA}{\mathcal{A}}

\newcommand{\cC}{\mathcal{C}}
\newcommand{\cD}{\mathcal{D}}
\newcommand{\cE}{\mathcal{E}}
\newcommand{\cF}{\mathcal{F}}
\newcommand{\cG}{\mathcal{G}}

\newcommand{\cM}{\mathcal{M}}
\newcommand{\cN}{\mathcal{N}}
\newcommand{\cP}{\mathcal{P}}
\newcommand{\cR}{\mathcal{R}}

\newcommand{\cX}{\mathcal{X}}

\newcommand{\sfC}{\mathsf{C}}
\newcommand{\sfE}{\mathsf{E}}

\newcommand{\sS}{\mathscr{S}}

\DeclareMathOperator{\Real}{Re}
\DeclareMathOperator{\Imaginary}{Im}
\DeclareMathOperator{\Span}{Span}

\theoremstyle{plain}
\newtheorem{thm}{Theorem}
\newtheorem{lem}{Lemma}
\newtheorem*{prop}{Proposition}

\theoremstyle{definition}
\newtheorem*{rem}{Remark}
\newtheorem*{defn}{Definition}

\begin{document}

\title{A dimensionally continued Poisson summation formula}

\author{Nathan~K.~Johnson-McDaniel}

\address{Institute for Gravitation and the Cosmos\\
  Center for Particle and Gravitational Astrophysics\\
  Department of Physics\\
  The Pennsylvania State University\\
  University Park, PA 16802, USA} 

\address{Theoretisch-Physikalisches Institut\\
  Friedrich-Schiller-Universit{\"a}t\\
  Max-Wien-Platz 1\\
  07743 Jena, Germany}

\email{nathan-kieran.johnson-mcdaniel@uni-jena.de}

\begin{abstract}
We generalize the standard Poisson summation formula for lattices so that it
operates on the level of theta series, allowing us to introduce
noninteger
dimension parameters (using the dimensionally continued Fourier
transform).
When combined with one of the proofs of the Jacobi
imaginary transformation of theta functions that does not use the
Poisson summation formula, our proof of this generalized Poisson
summation
formula also provides a new proof of the standard Poisson summation
formula for dimensions
greater than $2$
(with appropriate hypotheses on the
function being summed). 
 In general, our methods work to establish the (Voronoi)
summation formulae associated with functions satisfying 
(modular) transformations of the Jacobi imaginary type by means of a density argument
(as opposed to the usual Mellin transform approach). In particular, we construct a family of
generalized theta series from Jacobi theta functions from which these summation formulae can be obtained. This family
contains several families of modular forms, but is significantly more general than any of them. Our result also relaxes several
of the hypotheses in the standard statements of these summation formulae. The density result we prove for Gaussians in the Schwartz space may be
of independent interest.
\end{abstract}



\maketitle

\section{Introduction}

Consider a lattice $\Lambda\subset\R^n$ and a sufficiently well-behaved
function $F:\R^n\to\R$. [Taking $F$ to belong to the
Schwartz space $\sS(\R^n)$ is sufficient, and is what we shall do in our later
generalization.] The standard Poisson summation formula then says that
\<\label{nPSF}
\sum_{k\in\Lambda}F(k) = \frac{1}{\sqrt{\det\Lambda}}\sum_{p\in\Lambda^*}\tilde{F}(p).
\?
Here $\Lambda^*$ is the lattice dual to $\Lambda$, $\det\Lambda$ denotes the volume of a Voronoi
cell of $\Lambda$, and
\<\nonumber
\tilde{F}(p):=\int_{\R^n}F(x)e^{-2\pi i x\cdot p}dx
\?
denotes the Fourier transform of $F$. (We use a tilde here so that we can
reserve the circumflex for our more general, dimensionally continued Fourier
transform.) We wish to construct a dimensionally
continued version of this result.

This problem was originally inspired by
a condensed matter physics investigation involving the dimensional continuation
of electrostatic lattice sums, computed using the Ewald method (see,
e.g.,~\cite{JR} for a modern exposition of this method), as described in~\cite{J-MO1}. However, the ensuing
discussion is a purely mathematical offshoot of this investigation. For one
thing, the results we were able to prove do not include the physically
relevant case of a slowly decaying function, even though we have numerical
evidence that the results still hold in this case. Nevertheless, the methods used here might
also be applicable to the dimensional regularization of lattice sums:
See~\cite{CG} for an approach using zeta functions and the Mellin transform.

In addition to its use in performing lattice sums, the Poisson summation formula is a prominent
tool in most other problems involving lattices, from point counting problems in number theory
(discussed by, e.g., Miller and Schmid~\cite{MS}), to the mathematical theory of diffraction---see, e.g., the review by Lagarias~\cite{Lagarias}. While we do not consider diffraction in this
article, it is possible that our results might have some application to diffraction problems, or that
results in the theory of diffraction might provide some further avenues for generalization of the
results presented here. In particular, the quasicrystal summation formula given in Theorem~2.9
of Lagarias~\cite{Lagarias} seems to offer attractive possibilities for generalization. Additionally, Baake, Frettl{\"o}h, and Grimm~\cite{BFG} apply an integer dimension result that is very similar to
our formula to a problem in diffraction.


If we specialize to the case where $F$ is a radial function, we
can obtain a dimensionally continued Poisson summation formula in terms of
the lattice's theta series---see Theorem~\ref{Zd-thm} for the (particularly
simple) version for $\Z^d$ ($d \geq 1$). One can then generalize this result to a reasonably large
family of theta series-like functions constructed out of linear combinations of products of powers of
Jacobi theta functions; the generalization is given in Theorem~\ref{General_theorem}. This family
of generalized theta series includes all the theta series of lattices given in Chap.~4 of
Conway and Sloane~\cite{CS} (except for the general forms of those for
the root lattice $A_d$ and its translates), as well as several families of modular forms (including
all the classical modular forms of even weight), as discussed by
Rankin~\cite{Rankin}, and more compactly stated as Theorem~12 in Pache~\cite{Pache}. However, the family we consider is considerably more general, since we only
require integer powers to reproduce the aforementioned theta series from Conway and Sloane
and the families of modular forms from Rankin, while here the powers can be arbitrary
nonnegative real numbers.


This relation between modular transformations and
summation formulae is not new. In fact,
much of the basic result we obtain was first given by Bochner~\cite{Bochner}.  [The summation formula is stated
in Bochner's Eq.~(80) in an integral form and in his Eq.~(129) using sums; a nice form of sufficient hypotheses on the
function being summed are given in Bochner's Theorem~9. Also see Theorem~4.9.1 in~\cite{Bochner_HA} for a more streamlined presentation of the summation formula (using sums), though with some
obvious typos, and without a
statement of the hypotheses concerning the function being summed.] Our method of proof bears some
similarities to Bochner's, though Bochner uses Fourier-Laplace transforms as the basic analytic
tool, while we use a density result. (However, Bochner's method allows for weaker hypotheses.)
Our overall focus is also somewhat different, in that we are concerned with the summation
formulae derived from a specific family of generalized theta series, while Bochner only
considers a general modular relation. As far as we know, this specific construction of a family of
dimensionally continued Poisson summation formulae has not been considered previously.
While Guinand gives an
explicit analogue to our Theorem~\ref{Zd-thm} in Sec.~10.5 of~\cite{Guinand_AM}, this is
only stated for integer dimensions, even though Guinand cites Bochner as a reference for the
formula. (Guinand also cites his own general summation formula from~\cite{Guinand_PCPS}, which is only given for integer dimensions.) Note that we only became aware of the work of Bochner and Guinand
in the process of publication, through the reference to~\cite{Guinand_AM} in the review by
Lagarias~\cite{Lagarias} that the referee suggested to us. 

In more recent work, the relation has reached its most
refined form in the association between automorphic forms and Voronoi summation formulae (see,
e.g., Miller and Schmid~\cite{MS} for a review of recent work). However, work on this correspondence first arose in 
the context of transformations related to the functional equation for the zeta function, inspired
by a question from Voronoi on analogues of the Poisson summation formula---see, e.g.,~\cite{Sklar}
and references therein.
(We call particular attention to the work of Ferrar~\cite{Ferrar37}; see Theorem 10.2.17 in Cohen~\cite{Cohen} for a more modern discussion of a very similar result.)
Additionally, Baake,
Frettl\"{o}h, and Grimm~\cite{BFG} give a (distributional) radial Poisson
summation formula in their Theorem~3 in a form that is very similar to our
dimensionally continued form. However, they do not show how to dimensionally
continue the
lattice (or, indeed, mention theta functions explicitly), and their proof (which
relies on the standard Poisson summation formula) only holds for integer
dimensions.
There are also discussions of similar formulae---these derived from modular
transformations---at the beginning of Chap.~4 of Iwaniec and
Kowalski~\cite{IK}, and in Sec.~10.2 of Huxley~\cite{Huxley}---what Huxley
terms the Wilton summation formula. These formulae are presented in what appears to be
a dimensionally continued form, though their hypotheses assume integer dimensions. Regardless,
the summation formula in Iwaniec and Kowalski and Huxley's Wilton summation
formula are
derived from cusp forms, while our result (in the language of modular forms) does not require the
vanishing of the constant term in the form's Fourier series (the defining characteristic of a cusp form).



We 
make considerable use of various special functions in the following discussion, so
we provide references for the conventions we use. The dimensionally
continued theta series for which our summation formula holds will all be constructed out of Jacobi 
theta functions. An appropriate introduction to these functions for our purposes is given in Chap.~4
of Conway and Sloane~\cite{CS} [though there is an important difference in our notation, as
discussed below Eq.~\eqref{thetas}]. We shall also encounter the confluent
hypergeometric limit function and the Hermite functions (plus a brief appearance by a Bessel
function). These are discussed in an appropriate way by Andrews, Askey, and Roy~\cite{AAR}
(though we slightly streamline their general hypergeometric function notation, since we only
consider the confluent hypergeometric limit function). 
In general, the Wolfram Functions Site~\cite{WFS} is a good resource for
information about the various special functions we employ.

Additionally, we shall need the apparatus of tempered distributions introduced by
L.~Schwartz~\cite{Schwartz}. A convenient overview of
the properties we need is given in Sec.~V.3 of Reed and Simon~\cite{RS}. (The Appendix to that
section also provides another exposition of the $N$-representation we use.) 

\section{Ingredients}

\subsection{Theta series}
\label{ss:ts}

Here we recall various facts about theta series and theta functions
that we shall need for the rest of our discussion, following Chap.~4 of
Conway and Sloane~\cite{CS}. The theta series of a lattice $\Lambda$ is defined
by
\[
\Theta_\Lambda(q) := \sum_{k\in\Lambda}q^{|k|^2}.
\]
(This is often treated as a formal series, but converges for $q\in\C$, $|q|<1$.)
The utility of the theta series stems from the fact that the coefficient of
$q^l$ in the expansion of
$\Theta_\Lambda(q)$ in powers of $q$ gives the number of points in the intersection of the
lattice and a sphere of radius $\sqrt{l}$ centred at the origin. Thus, if we write
\<\label{Theta}
\Theta_\Lambda(q) =: \sum_{l=0}^\infty N_l q^{A_l},
\?
then
\[
\sum_{k\in\Lambda}f(|k|) =
\sum_{l=0}^\infty N_l f(\sqrt{A_l}).
\]
[\emph{Nota bene}: We have written the radial function $F$ as $f(|\cdot|)$, and
shall only consider these radial parts in the sequel.]
Examples of dimensionally continued theta series for families of lattices include the
$d$-dimensional cubic lattice $\Z^d$,
with $\Theta_{\Z^d}(q) = \vartheta_3^d(q)$, and the root lattice $D^d$, with $\Theta_{D^d}(q) =
[\vartheta_3^d(q) + \vartheta_4^d(q)]/2$. (See Chap.~4 in Conway and Sloane~\cite{CS} for further
examples.) Here
\<\label{thetas}
\vartheta_2(q) := 2q^{1/4}\sum_{l=1}^\infty q^{l^2 - l}, \qquad \vartheta_3(q) := 1+2\sum_{l=1}^\infty q^{l^2}, \qquad \vartheta_4(q) := 1 + 2\sum_{l=1}^\infty (-q)^{l^2}
\?
are Jacobi theta functions (where $\vartheta_2$ is defined for future use).

\emph{Nota bene}: It is often
customary to take theta functions and theta series to be functions of a complex variable $z$,
instead of the nome $q=e^{i\pi z}$ that we have used here.
We have chosen to regard the nome as
fundamental since we are primarily interested in the expansions of these functions in powers of
$q$. However, when discussing transformations of these functions, it is considerably
more convenient to regard them as functions of $z$. On the few occasions where we do this, we
shall use an overbar to denote the difference, e.g., $\bar{\Theta}_\Lambda(z) :=
\Theta_\Lambda(e^{i \pi z})$. (In the literature on summation formulae
derived from automorphic forms, one thinks of our expansions of $\Theta_\Lambda$ in $q$ as the
Fourier coefficients of $\bar{\Theta}_\Lambda$.)

Since the Poisson summation formula involves the dual lattice, we need to know how to obtain its
theta series. This is given by the Jacobi formula [Eq.~(19) in Chap.~4 of Conway and
Sloane~\cite{CS}], which states that
\<\label{JTF1}
\bar{\Theta}_{\Lambda^*}(z) =
\sqrt{\det\Lambda}(i/z)^{d/2}\bar{\Theta}_\Lambda(-1/z),
\?
where $d$ is the dimension of the lattice.
The Jacobi formula is typically proved using the Poisson summation formula [see, e.g., the
discussion leading up to our Eq.~\eqref{JTF2}]. However, all we need in our
discussion is the
intimately related Jacobi imaginary transformation of the Jacobi theta 
functions
(also known as the modular identity or reciprocity formula for the theta functions), i.e.,
\<\label{JIT}
\bar{\vartheta}_2(-1/z) = (z/i)^{1/2}\bar{\vartheta}_4(z), \qquad\qquad
\bar{\vartheta}_3(-1/z) = (z/i)^{1/2}\bar{\vartheta}_3(z).
\?
(The first of these is also true with the labels $2$ and $4$ switched.)
The standard proof of these identities is a direct application of the Poisson
summation formula, but there are alternative proofs that are independent of it.
For instance, one such proof is given in Sec.~21.51 of
Whittaker and Watson~\cite{WW}, while Bellman's text~\cite{Bellman}
discusses several others---see, in particular, Sec.~30 for Polya's derivation---in addition to the
standard Poisson summation version (in Sec.~9). Our discussion will
thus be independent of the
standard Poisson summation formula (with the exception of a brief appeal to establish
Theorem~\ref{Zd-thm} for $d = 1$).

\subsection{The dimensionally continued Fourier transform}
\label{d-FT}

We also need to dimensionally continue the Fourier transform. Stein and Weiss give a dimensionally
continued version of the Fourier transform for radial functions in Theorem~3.3 of Chap.~IV of~\cite{SW}, viz.,
\<\label{f-hat}
\begin{split}
\hat{f}(p) &:= 2\pi p^{-(d-2)/2}\int_0^\infty f(r)J_{(d-2)/2}(2\pi p r)r^{d/2}dr\\
&= \frac{2\pi^{d/2}}{\Gamma(d/2)}\int_0^\infty
f(r){}_0F_1(d/2;-\pi^2p^2r^2)r^{d-1}dr.
\end{split}
\?
(This reduces to the standard Fourier transform for a radial function when $d\in\N$.)
Here the first equality gives the expression from Stein and Weiss ($J_k$ is a Bessel function) and the second
gives an equivalent (perhaps slightly neater) expression in terms of the confluent hypergeometric limit function~${}_0F_1$.
The hypergeometric expression has the advantage of only involving one
appearance of $p$ (and being manifestly regular at $p = 0$ for all $d \ge 1$),
in addition to showing the
$d$-dimensional
polar coordinate measure for radial functions explicitly. We shall thus use
the hypergeometric expression exclusively in the sequel. (One can obtain the hypergeometric expression using the Stein and
Weiss derivation---the only difference is that one uses a different special function to evaluate the final
integral.\footnote{The integral representation for ${}_0F_1$ we used is
07.17.07.0004.01 on the Wolfram Functions Site~\cite{WFS}.})

For this expression to be well-defined, it is sufficient
to take $d \ge 1$: One assumes $d>1$ when using integral representations to express the result in terms of
either of the two given special functions, and can also check that the integral is convergent for all $p\in\R$
in that case, provided that $f\in L^1(\R_+)$. Additionally, Eq.~\eqref{f-hat} reduces to the expected
expression for $d=1$ [using $J_{-1/2}(z) = \sqrt{2/\pi z}\cos z$ or ${}_0F_1(1/2;-z) = \cos(2\sqrt{z})$].\footnote{These identities are 03.01.03.0004.01 and
07.17.03.0037.01, respectively, on the Wolfram Functions Site~\cite{WFS}.}
(Stein and Weiss restrict to $d\ge 2$ so that the integral they use in their derivation is well-defined, since they are only considering integer dimensions.)
This restriction to $d > 1$ is necessary for other parts of our
discussion, though we have numerical evidence that it can be relaxed.

The following result is central to understanding why this dimensionally
continued Fourier transform agrees with the dimensional continuation of the theta series. 

\begin{lem}\label{G-FT}
For $d \ge 1$, the dimensionally continued Fourier transform (for radial functions) defined in
Eq.~\eqref{f-hat} takes a Gaussian
$\cG_\alpha(r) := e^{-\alpha r^2}$, $\Real \alpha > 0$ to another
Gaussian, $\hat{\cG}_\alpha(p) = (\pi/\alpha)^{d/2}e^{-\pi^2p^2/\alpha}$.
\end{lem}

\begin{rem}
Intuitively, this result follows from dimensionally continuing the well-known
integer dimension result. We should get the same result from direct
calculation using Eq.~\eqref{f-hat} since that expression was obtained using
the same dimensional continuation procedure.
\end{rem}

\begin{proof}
The case $d = 1$ is classical. For $d > 1$, we use ${}_0F_1$'s defining series,
\<\label{0F1-series}
{}_0F_1(d/2;-\pi^2p^2r^2) = \sum_{n=0}^\infty\frac{(-\pi^2p^2r^2)^n}{(d/2)_n n!}
\?
[$(\cdot)_n$ denotes the Pochhammer symbol] and integrate term-by-term,
evaluating each integral using the gamma function.\footnote{The Maclaurin
series for ${}_0F_1$ is 07.17.02.0001.01 on the
Wolfram Functions Site~\cite{WFS}.} The resulting series
is the Maclaurin series for the expression we gave for $\hat{\cG}_\alpha$. The
term-by-term integration is justified by the Lebesgue dominated convergence
theorem. To see this, we use the same integral representation
for ${}_0F_1$ used in the derivation of Eq.~\eqref{f-hat}, which gives, for any
$N\in\N$,
\<\nonumber
\left|\sum_{n=0}^N\frac{(-\pi^2p^2r^2)^n}{(d/2)_n n!}\right| \le
{}_0F_1(d/2;\pi^2p^2r^2)
\le K\cosh(2\pi pr)\int_0^1(1-t^2)^{(d-3)/2}dt,
\?
where $K > 0$ is a constant.\footnote{\emph{Nota
bene}: We denote the set of positive integers by $\N$, and the set of
nonnegative integers by $\N_0$.} This allows us to apply the dominated
convergence theorem, since the integral in the final term is finite for
$d>1$ and $\int_0^\infty \cosh(2\pi pr)|e^{-\alpha r^2}|dr$ is finite for
$\Real\alpha>0$.
\end{proof}

\begin{rem}
The importance of this result to our discussion comes in its use in 
obtaining the integer dimension Jacobi transformation formula (and thus also
the Jacobi imaginary transformations of the Jacobi theta functions) via the
standard Poisson summation
formula: For a lattice $\Lambda$ of dimension $n\in\N$, we have (taking
$\Imaginary z > 0$ so that everything converges)
\<\nonumber
\bar{\Theta}_\Lambda(z) := \sum_{k\in\Lambda}e^{i\pi z|k|^2} =
\frac{1}{\sqrt{\det\Lambda}}\left(\frac{i}{z}\right)^{n/2}\sum_{p\in\Lambda^*}e^{-i\pi|p|^2/z}
= \frac{1}{\sqrt{\det\Lambda}}\left(\frac{i}{z}\right)^{n/2}\bar{\Theta}_{\Lambda^*}(-1/z), 
\?
which can be written as
\<\label{JTF2}
\bar{\Theta}_{\Lambda^*}(z) = \sqrt{\det\Lambda}(i/z)^{n/2}\bar{\Theta}_\Lambda(-1/z),\?
the Jacobi
transformation formula. We thus expect that the dimensionally continued dual theta series that we obtain
using this formula will agree with the dimensionally continued Fourier transform to give a dimensionally
continued Poisson summation formula.
\end{rem}

\section{The dimensionally continued Poisson summation formula for $\Z^d$}

With these results in hand, we can thus write Eq.~\eqref{nPSF} [for a radial
function $F =: f(|\cdot|)$] as
\<\nonumber
\sum_{l=0}^\infty N_l f(\sqrt{A_l}) =
\frac{1}{\sqrt{\det\Lambda}}\sum_{l=0}^\infty N^*_l \hat{f}(\sqrt{A^*_l}),
\?
where the starred quantities come from writing the theta series of
$\Lambda^*$ in the power series form given by Eq.~\eqref{Theta},
and we calculate $\hat{f}$ by taking the dimension parameter $d$ to be the dimension of the lattice. (As we shall see
later, what is important is that the $d$ one uses here is the same $d$
that appears in the Jacobi
transformation formula.) It is clear that this equality holds when $d\in\N$, by the standard
Poisson summation formula. What is perhaps surprising is that the equality still holds
for, e.g.,
$\Lambda = \Z^d$, with $d\in\R$ ($d \geq 1$). We shall first prove the result for this
simple case ($\Z^d$ is self-dual, $\det\Z^d=1$, and $A_l = l$), where it becomes

\begin{thm}\label{Zd-thm}
If $f\in\sS^\sfE(\R)$ (i.e., $f$ is an even Schwartz function) and $d\ge 1$, then
\<\label{ZPSF}
\sum_{l=0}^\infty
N_lf(\sqrt{l}) = \sum_{l=0}^\infty
N_l\hat{f}(\sqrt{l}),
\?
where the $N_l$ are given by the power series expansion of the theta series of $\Z^d$, viz.,
\<\label{Theta-Zd}
\Theta(q) =
\vartheta_3^d(q) = \left[1+2\sum_{k=1}^\infty q^{k^2}\right]^d
=: \sum_{l=0}^\infty N_lq^l,
\?
and $\hat{f}$ is computed using Eq.~\eqref{f-hat}.
\end{thm}
However, the simplifications are primarily notational. As we shall see in the
discussion in Sec.~\ref{sec:generalization}, the proof
works with minimal modifications for a much larger class of $\Theta$s,
including functions that cannot be the theta series
of a lattice (even though they have an integer dimension parameter).

\begin{rem}
The restriction that $f$ be an even function should not be surprising: In
integer dimensions, it corresponds to the lack of a cusp at the origin for the
full radial function $F = f(|\cdot|)$. Moreover, as Miller and Schmid note~\cite{MS}, the
standard one-dimensional Poisson summation formula is a trivial $0 = 0$ for odd
functions.
\end{rem}

\section{A Schwartz space density result}

Since the proof proceeds by noting that the desired formula holds almost
trivially for the Gaussians from Lemma~\ref{G-FT}, and then extends to an
interesting set of functions [viz., $\sS^\sfE(\R)$]
by density, we start by establishing the requisite density result.

\begin{lem}\label{density}
$\overline{\Span\{x\mapsto e^{-\alpha x^2}|\alpha>0\}}^{\sS(\R)} =
\sS^\sfE(\R)$ [i.e., the Schwartz space closure of the given family of
Gaussians is all the even
Schwartz functions].
\end{lem}

\begin{proof}
We shall prove this by showing that
\<\nonumber
\cX := \Span\{x\mapsto e^{-\alpha x^2}|\alpha>0\} +
\Span\{x\mapsto xe^{-\alpha x^2}|\alpha>0\}
\?
is dense in $\sS(\R)$, so its even part,
$\Span\{x\mapsto e^{-\alpha x^2}|\alpha>0\}$, is thus dense in
$\sS^\sfE(\R)$. We shall use Corollary
IV.3.14 from Conway's text~\cite{Conway}, which states that a linear manifold (here $\cX$) is dense in a locally
convex topological vector space [here $\sS(\R)$] if and only if the only element of the dual of the
topological vector space that vanishes on all elements of the linear
manifold is the zero element.

It is most convenient to proceed by
identifying
$\sS(\R)$ with a sequence space, following Simon~\cite{Simon}. (There is an alternative
presentation of these results in the Appendix to Sec.~V of Reed and Simon~\cite{RS}.) Here the
sequence space is given by the coefficients of the
Hermite function expansion of elements of $\sS(\R)$, and provides a
particularly
nice characterization of the tempered distributions [the elements of
$\sS'(\R)$, the dual of $\sS(\R)$]. Namely, if $a_n$ are the Hermite
coefficients of $f\in\sS(\R)$ [i.e., $a_n := \int_\R f(x)h_n(x)dx$, where $h_n$ is
the $n$th Hermite function], then $\varphi\in\sS'(\R)$ can be written as
$\varphi(f) = \sum_{n=0}^\infty c_na_n$, where $c_n$ are the Hermite
coefficients of $\varphi$, with $|c_n| \leq C(1+n)^m$ for some
$C, m > 0$. (This is Theorem~3 in Simon~\cite{Simon}.) Note that Simon defines the Hermite functions
to be $L^2$ normalized,
so, we have, from the first equation in Sec.~2 of Simon,\footnote{\emph{Nota bene}: Simon defines the $h_n$
without the factor of $(-1)^n$ (that here comes from our $H_n$). We have included the $(-1)^n$
for notational simplicity (since we use the standard convention for the Hermite
polynomials). This does not have any effect on Simon's Theorem~3, since it
simply amounts to a sign change of the odd Hermite coefficients.}
\[
h_n(x) := \frac{e^{-x^2/2}}{\sqrt{\pi^{1/2}2^n n!}}H_n(x), \qquad\qquad
H_n(x) := (-1)^ne^{x^2}\frac{d^n}{dx^n}e^{-x^2},
\]
where the $H_n$ are the Hermite polynomials, with generating
function\footnote{This is 05.01.11.0001.01 on the Wolfram Functions
Site~\cite{WFS}.}
\[
\sum_{n=0}^\infty H_n(x)\frac{t^n}{n!} = e^{2tx - t^2}.
\]

We can now use this generating function to show that
the Hermite coefficients of $x \mapsto e^{-\alpha x^2}$ are
given by
\<\nonumber
a_n = \cN_n\frac{d^n}{dt^n}\left.\left[\int_\R e^{-(\alpha x^2 + x^2/2
    - 2tx + t^2)}dx\right]\right|_{t = 0} =
\cN_n\frac{d^n}{dt^n}\left.[\sqrt{\pi\beta}e^{(\beta
      - 1)t^2}]\right|_{t = 0},
\? 
where $\cN_n := (\pi^{1/2}2^n n!)^{-1/2}$ is the Hermite functions' normalization factor and
$\beta := 1/(\alpha + 1/2)$. We thus have $a_{2n} = \cN_{2n}(\pi/\beta)^{1/2}(\beta - 1)^n/n!$, $a_{2n+1} = 0$,
by the series expansion of the exponential. [We used Lemma~2.2 in Chap.~VIII of
Lang~\cite{Lang} to interchange differentiation and integration. We only need to consider
the case where $t$ lies in some neighbourhood of $0$, so the $t$-derivatives of the integrand are each bounded by
a polynomial in $x$ times a Gaussian in $x$ (for all $t$ in the neighbourhood), and
those functions of $x$ are integrable over $\R$.]
Similarly, the Hermite coefficients of $x \mapsto
xe^{-\alpha x^2}$ are $b_{2n} = 0$ and $b_{2n+1} =
\cN_{2n+1}(\pi/\beta^3)^{1/2}(\beta - 1)^n/n!$. Thus, we consider
\[
E_{\beta,\pm}(x) := (\beta/\pi)^{1/2}e^{-\alpha x^2}
\pm (\beta^3/\pi)^{1/2}xe^{-\alpha x^2},
\]
which has Hermite
coefficients of $(\pm 1)^n\cN_n(\beta - 1)^{\lfloor n/2\rfloor}/\lfloor
n/2\rfloor!$, where $\lfloor\cdot\rfloor$ denotes the greatest integer
less than or equal to its argument.

Now, for any $\varphi\in\sS'(\R)$, $\cE_\pm(\beta) := \varphi(E_{\beta,\pm})$ is a holomorphic function of $\beta$.
To see this, we note that
\<
\cE_\pm(\beta) = \sum_{n=0}^\infty(\pm 1)^nc_n\cN_n\frac{(\beta - 1)^{\lfloor n/2\rfloor}}{\lfloor
n/2\rfloor!} = \sum_{n=0}^\infty(\cN_{2n}c_{2n} \pm \cN_{2n+1}c_{2n+1})\frac{(\beta - 1)^n}{n!},
\?
where $c_n$ are the Hermite coefficients of $\varphi$. Since the $c_n$
are bounded by a polynomial in $n$, the series converges for all $\beta\in\C$,
giving holomorphy.
Thus, if $\cE_\pm(\beta) = 0$ for all $\beta$ in an interval
(as is the case here), then all of $\cE_\pm$'s power series
coefficients are zero. Applying this result to the two choices of sign, we
obtain (since the $\cN_n$ are never zero) $c_n = 0$ $\forall$ $n\in\N_0$ $\Rightarrow$ $\varphi
\equiv 0$, which thus proves the lemma.
\end{proof}

\begin{rem}
This result may be of wider
applicability, particularly in harmonic analysis, due to the ubiquity of the
Gaussian. We thus note that the proof of the lemma shows that $\alpha$ need merely
belong to some subset of the right half-plane with an accumulation point to
guarantee density.
One could have also proved this result more abstractly (and without recourse to
the Hermite expansion) by a slightly indirect application of the
Stone-Weierstrass theorem, though the basic Hahn-Banach argument (contained in
the Corollary from Conway we use) remains the
same.\footnote{Personal communication from John Roe.}
\end{rem}

\section{Proof of Theorem~\ref{Zd-thm}}

We first note that Eq.~\eqref{ZPSF} is clearly true for $d=1$ (indeed, $d\in\N$) by the standard Poisson summation
formula for lattices (applied to $\Z^d$). To prove the result for
$d>1$, we shall first establish that it holds for the
Gaussians
from Lemma~\ref{G-FT}, and then show that the equality still holds in the limit
in the Schwartz space topology.
The control afforded by demanding convergence in the Schwartz space
makes this quite straightforward. The primary result that needs to be shown is that two functions that are $\epsilon$-close in the Schwartz space topology
have dimensionally continued Fourier transforms that are $C\epsilon$-close in
a given Schwartz space seminorm (where the constant $C$ depends on the seminorm
under consideration, as well as $d$).

To show that Eq.~\eqref{ZPSF} holds when $f = \cG_\alpha$, we first consider the left-hand side and note that
\<\label{Taylor}
\sum_{l=0}^\infty N_le^{-\alpha l}=
\Theta(e^{-\alpha}).
\?
Convergence is guaranteed because $\Theta$ is analytic inside the unit
disk.
[To see that $\Theta$ is analytic inside the
unit disk, note that $\vartheta_3$ is analytic there, and, moreover, nonzero, so its $d$th power is
analytic, as well. It is easiest to see that $\vartheta_3$ is nonzero inside the unit disk from its
infinite product expansion, given in, e.g., Eq.~(35) in Chap.~4 of Conway and Sloane~\cite{CS}.] Using Lemma~\ref{G-FT},
the right-hand side of Eq.~\eqref{ZPSF} becomes
\<\nonumber
\left(\frac{\pi}{\alpha}\right)^{d/2}\sum_{l=0}^\infty N_l e^{-\pi^2l/\alpha} = \left(\frac{\pi}{\alpha}\right)^{d/2}\Theta(e^{-\pi^2/\alpha}).
\?
Now, the Jacobi imaginary transformation for $\vartheta_3$ [Eq.~\eqref{JIT}] implies that $(\pi/\alpha)^{d/2}\Theta(e^{-\pi^2/\alpha}) = \Theta(e^{-\alpha})$, so we have thus established
the result for $\cG_\alpha$.

We shall now show that this equality continues to hold in the limit. The equality is clearly true
for any finite linear combination of the Gaussians from Lemma~\ref{G-FT}, so we use Lemma~\ref{density}
to approximate an arbitrary $f\in\sS^\sfE(\R)$ by a
finite
linear combination of these Gaussians, $g$. Specifically, we have
$\|f - g\|_{n,m} < \epsilon$ $\forall$ $n,m\in\N_0$, where $\|f\|_{n,m} := \sup_{x\in\R}|x^nf^{(m)}(x)|$
is the family of seminorms that gives the Schwartz space topology. (We denote the $m$th derivative
of $f$ by $f^{(m)}$.) We wish to bound the difference
between the two sides of Eq.~\eqref{ZPSF} by a constant times $\epsilon$.
We have
\<\label{bound}
\left|\sum_{l=0}^\infty N_lf(\sqrt{l}) - \sum_{l=0}^\infty N_l\hat{f}(\sqrt{l})\right| \le
\left|\sum_{l=0}^\infty N_l(f-g)(\sqrt{l})\right| + \left|\sum_{l=0}^\infty N_l(\hat{f}-\hat{g})(\sqrt{l})\right|,
\?
where we used the fact that the dimensionally continued Poisson summation
formula holds for $g$, along with the triangle inequality.
We can bound the two sums on the right-hand side by constants times
$\epsilon$ using the assumption about the closeness of $f$
to $g$ in the Schwartz space topology and the fact that $N_l$ grows at most
polynomially with $l$. The latter fact also shows that the two sums on the left converge for $f\in\sS(\R)$.

\subsection{Bounds on the growth of $N_l$ and on the right-hand side of Eq.~\eqref{bound}}\label{Nl_bounds}

\begin{lem}
The $N_l$ defined in Eq.~\eqref{Theta-Zd} are polynomially bounded. Specifically, we have
\[
|N_l| \le 2^d(1+d/l)^l(1+l/d)^d \le C_dl^d,
\]
where $C_d>0$ is some constant (and the second inequality only holds for $l \geq 1$).
\end{lem}

\begin{proof}
Recalling that $\Theta$ is analytic inside the unit disk, we can apply Cauchy's integral formula to the contour $\sfC_R$, a
circle of radius $R\in(0,1)$, centered at the origin (and oriented
counterclockwise), to obtain
\<\nonumber
|N_l| = \left|\frac{1}{2\pi i}\int_{\sfC_R}\frac{\vartheta_3^d(z)}{z^{l+1}}dz\right| = \frac{1}{2\pi}\left|\int_0^{2\pi}\frac{\vartheta_3^d(Re^{i\theta})}{R^le^{il\theta}}d\theta\right| \le \frac{2^d}{R^l(1-R)^d}.
\?
Here we have used $|\vartheta_3(q)| \le 2/(1-|q|)$ (for $|q| < 1$, obtained
using the geometric series). The right-hand side attains its minimum
[for $R\in(0,1)$] at $R=l/(l+d)$, so we have the desired result.  To obtain the second inequality,
we use the fact that $(1 + 1/r)^r < e$ for $r>0$.
\end{proof}

\begin{rem}
While this bound is easy to
obtain and is all that
is necessary for our purposes, it is by no means optimal. For instance, for integer $d$, we can apply the
Hecke bound for modular forms (e.g., Theorem~4.5.3 in Rankin~\cite{Rankin}) to conclude that
$N_l = O(l^{d/2})$. (Other of Rankin's results---stated as Theorem~12 in Pache~\cite{Pache}---show that the Hecke bound stated by Rankin is applicable to $\vartheta_3^d$ for $d\in\N$.)
\end{rem}

If we write $h := f - g$, then this bound implies that $|N_lh(\sqrt{l})| \le C_dl^d|h(\sqrt{l})| \le \epsilon C_d/l^2$ (for $l \geq 1$), where the second inequality
follows from the fact that $h$ is $\epsilon$-close to $0$ in the Schwartz space topology. [Explicitly, we have
$|x^{2d + 4}h(x)| \le \epsilon $ $\forall$ $x > 1$ $\Rightarrow l^d|h(\sqrt{l})| \le \epsilon/l^2$ $\forall$ $l\in\N$. The first inequality comes from
noticing that for any $\gamma \ge 0$,
we have $|x^\gamma h(x)| \le
|x^{\lceil\gamma\rceil}h(x)| \le \epsilon$ for $x \geq 1$, where $\lceil\cdot\rceil$ denotes
the smallest integer greater than or equal to its argument.]
We shall show that
$|p^{2n}\hat{h}(p)| \le K_d\epsilon$ $\forall$ $n\in\N$, $p\in\R$ (where $K_d$ is some $n$-dependent constant), so we have $|p^{2d + 4}\hat{h}(p)| \le K_d\epsilon$ $\forall$ $p\in\R$.
We can thus apply the same argument to
the second sum and hence bound both sums by constants times $\epsilon$
(since $\sum_{l=1}^\infty l^{-2}$ is finite), showing that the dimensionally
continued Poisson summation
formula is true in the limit [since we will have shown that the right-hand side of Eq.~\eqref{bound} is bounded by a constant times $\epsilon$].

\subsection{Bound on $|p^{2n}\hat{h}(p)|$}

To prove the bound on $|p^{2n}\hat{h}(p)|$, we first
dimensionally continue some standard Fourier results.

\begin{lem}\label{eig-lem}
If we define the $d$-dimensional Laplacian for radial functions by 
\<\label{d-Laplacian}
\triangle_df(r) := f''(r) + \frac{d-1}{r}f'(r),
\?
then, for $d > 1$,
\begin{enumerate}[i)]
\item $\cF_p(r) := {}_0F_1(d/2;-\pi^2p^2r^2)$ satisfies
$\triangle_d\cF_p = -4\pi^2p^2\cF_p$, so
\item $\widehat{\triangle_d^nf}(p) =
(-1)^n(2\pi p)^{2n}\hat{f}(p)$ for $f\in\sS(\R)$.
\end{enumerate}
\end{lem}

\begin{proof}
Part \emph{i} follows from the fact that $y_a(r) := {}_0F_1(a;r)$ satisfies
$ry_a''(r) + ay_a'(r) = y_a(r)$.\footnote{This differential equation for
${}_0F_1$ is 07.17.13.0003.01 on the
Wolfram Functions Site~\cite{WFS}.} [Alternatively, it can be obtained by direct
calculation using
Eq.~\eqref{0F1-series}, justifying term-by-term differentiation using analyticity.] Part \emph{ii} is then obtained by induction, applying
Eq.~\eqref{f-hat} to $\triangle_d^{n-1}f$ and integrating by parts twice. [The
boundary terms at infinity vanish because $f\in\sS(\R)$; those at $0$ vanish
because $d>1$ (or cancel amongst themselves).]
\end{proof}

We can thus write $|p^{2n}\hat{h}(p)| = (2\pi)^{-2n}|\widehat{\triangle_d^nh}(p)|$. Then, since we shall show below that
$|r^k\triangle_d^nh(r)|\le\cD\epsilon$, where $\cD$
is some ($n$- and $d$-dependent constant), we obtain [using
Eq.~\eqref{f-hat} and the fact that ${}_0F_1(a;r)$ is a bounded
function of $r$, as was seen in the proof of Lemma~\ref{G-FT}]
\<\nonumber
\begin{split}
|p^{2n}\hat{h}(p)| &\le \cC\int_0^\infty |\triangle_d^nh(r)|r^{d-1}dr\\
&\le \cC\left[\int_0^1 |\triangle_d^nh(r)|dr + \int_1^\infty |\triangle_d^nh(r)|r^{d-1}dr\right]\\
&\le \cC\cD\left[1 + \int_1^\infty r^{d-1-s}dr\right]\epsilon,
\end{split}
\?
where  $\cC > 0$ is some ($n$- and $d$-dependent) constant and we used
$|r^k\triangle_d^nh(r)|\le\cD\epsilon$ with $k = 0$ and $k = s$. We can
choose $s > d$, so the integral in the final term is finite, thus giving the desired
result.

\subsection{Bound on $|r^k\triangle_d^nh(r)|$}

To see that $|r^k\triangle_d^nh(r)|$ is bounded by some ($n$- and $d$-dependent) constant
(called $\cD$ above), we first note
that we can use induction to write
\<\label{tri-express}
\triangle_d^nh(r) = \sum_{j=1}^{2n}a_j\frac{h^{(j)}(r)}{r^{2n-j}}
\?
for some ($n$- and $d$-dependent) constants $a_j$ (and an arbitrary
differentiable $h$). Thus, for $|r|\ge 1$, we have $|r^k\triangle_d^nh(r)|\le
\epsilon\sum_{j=1}^{2n}a_j$. For $|r| < 1$, matters are considerably more subtle, and we have to rely on the
fact that $h$ is even to see that $\triangle_d^nh$ remains bounded at the origin. The
argument goes as follows: We write $h = \cP + \cR$, where
$\cP$ is $h$'s $(2n)$th-degree Maclaurin polynomial (necessarily even, since $h$ is)
and $\cR$ is the associated remainder. We then have
$|r^k\triangle_d^nh(r)| \le |r^k\triangle_d^n\cP(r)| +
|r^k\triangle_d^n\cR(r)|$. Since $\triangle_d^n$ maps even polynomials to even polynomials
[as can be seen from Eq.~\eqref{tri-express}], $|r^k\triangle_d^n\cP(r)|$ is
bounded by a
($k$-, $n$-, and $d$-dependent) constant times $\epsilon$ for $|r|\le
1$. [Since
the coefficients of $\cP$ are given by derivatives of $h$,
they are bounded by constants times $\epsilon$, by hypothesis.]
To deal with $|r^k\triangle_d^n\cR(r)|$, we first need to establish an identity for derivatives of $\cR$, viz., (for $j \leq 2n$)
\<\nonumber
\cR^{(j)}(r) = \frac{h^{(2n+1)}(\xi_j)}{(2n+1-j)!}r^{2n+1-j},
\?
for some $\xi_j \in(0,r)$. One obtains this by comparing the $j$th derivative of $h = \cP + \cR$
with the $(2n-j)$th order Maclaurin expansion (with
Lagrange remainder) of $h^{(j)}$. The polynomial pieces are the same, while the remainder
pieces give the two sides of the equality. Combining this identity with Eq.~\eqref{tri-express},
we obtain
\<\nonumber
\triangle_d^n\cR(r) = r\sum_{j=1}^{2n}b_jh^{(2n+1)}(\xi_j),
\?
where the $b_j$ are ($n$- and $d$-dependent) constants. This
shows that $|r^k\triangle_d^n\cR(r)|$ is bounded
by an ($n$- and $d$-dependent) constant times $\epsilon$ for $|r|\le
1$, so
$|r^k\triangle_d^nh(r)|$ is, as well, proving the desired result, and hence
the theorem. \qed

\begin{rem}
The restrictions on $f$ and $d$ in the statement of the theorem are surely not
optimal: Bochner's proof of the basic summation formula holds for $d\in\R$, $d > 0$, and $f$ of only finite differentiability (see Theorem~9 in~\cite{Bochner}), and
there is numerical evidence that the given
result holds for $d\in\C$, $\Real d > 0$ and significantly less smooth $f$ [e.g., $f(r) = e^{-|r|^3}$, which is not covered by Bochner's results]. (The evidence also extends to the generalization given in Theorem~\ref{General_theorem} and is provided by a {\sc{Mathematica}}
notebook,
available online.\footnote{The notebook is available at {\url{http://gravity.psu.edu/~nathanjm/Dim_cont_PSF_test.nb}}.})
While one could use a slightly larger function space than $\sS^\sfE(\R)$
without any change
to the proof---the proof does not need control over $\|f - g\|_{n,m}$ for all $n$
and $m$---we did not investigate this in any detail: The resulting function
space would still require a fair amount of differentiability (while we have
numerical evidence that the formula remains true for at least some functions
with a cusp at the origin), and faster decay than the standard Poisson
summation formula (see, e.g., Corollary~2.6 in Chap.~VII of Stein and Weiss~\cite{SW}).  Moreover, the closure of the family of Gaussians in this
less restrictive topology would almost surely be more recondite than
$\sS^\sfE(\R)$.
\end{rem}

\section{Generalization of Theorem~\ref{Zd-thm}}
\label{sec:generalization}

Since there are other families of lattices with dimensionally continued theta
series besides $\Z^d$ (e.g., the root lattice $D^d$ mentioned in
Sec.~\ref{ss:ts}),
it is reasonable to expect that Theorem~\ref{Zd-thm} can be
generalized by replacing $\Theta$ with some more general function $\Upsilon$, which we shall term a generalized theta series. It
is not clear how to construct the most general such $\Upsilon$.\footnote{But note that
Ryavec characterizes all admissible $\Upsilon$s
(under certain assumptions) for $d = 1$ in~\cite{Ryavec}. We also call
attention to the work of C\'{o}rdoba~\cite{Cordoba, Cordoba_LMP}, who shows
that in integer dimensions, large classes of generalized Poisson summation
formulae arise from the standard Poisson summation formula applied to the finite disjoint union of (integer dimensional) lattices. (Note that Lagarias makes a slight correction to the statement of Theorem~2 of~\cite{Cordoba} in Theorem~3.7 of~\cite{Lagarias}.)} However, the template
provided by the theta series of other standard lattices (e.g., the ones given in Chap.~4
of Conway and Sloane~\cite{CS}) allows us to construct a reasonably general family of
generalized theta series out of finite linear combinations of products of the three Jacobi theta functions
($\vartheta_2$, $\vartheta_3$, and $\vartheta_4$) given in Eq.~\eqref{thetas}. This family will
contain all of the theta series of lattices and shifted lattices given in Chap.~4 of Conway and
Sloane~\cite{CS}, except for the general form of the theta series of the root lattice $A^d$ and
its translates.
In fact, theorems in Conway and
Sloane (Theorems~7, 15, and 17 in Chap.~7 and Theorem~5 in Chap.~8) show that the theta series of large classes of lattices can be written in such a form. 
However, the family of $\Upsilon$s is considerably more general,
since one only requires
$\lambda_m,\rho_m,\sigma_m\in\N_0$ to reproduce the theta series in Conway and Sloane, while here they can be arbitrary nonnegative real numbers. 
The $\Upsilon$s also contain several families of modular forms, as shown by Rankin~\cite{Rankin}, and concisely stated in Theorem~12 of Pache~\cite{Pache}.

Explicitly, we make the following

\begin{defn}
A \emph{generalized theta series} is a finite linear combinations of $\Upsilon_d$s of the form
\<
\Upsilon_d(q) := \prod_{m=1}^\cM\vartheta_2^{\lambda_m}(q^{s_m})\vartheta_3^{\rho_m}(q^{t_m})\vartheta_4^{\sigma_m}(q^{u_m}),
\?
with $\lambda_m,\rho_m,\sigma_m \geq 0$, $\sum_{m=1}^\cM(\lambda_m + \rho_m + \sigma_m) = d$, and $s_m,t_m,u_m\in\Q_+$.
\end{defn}

We thus have
\<\label{Upsilon_d_star}
\Upsilon_d^*(q) = \prod_{m=1}^\cM\frac{\vartheta_2^{\sigma_m}(q^{1/u_m})\vartheta_3^{\rho_m}(q^{1/t_m})\vartheta_4^{\lambda_m}(q^{1/s_m})}{\sqrt{s_m^{\lambda_m}t_m^{\rho_m}u_m^{\sigma_m}}},
\?
which we compute using
\[
\bar{\Upsilon}^*(z) := (i/z)^{d/2}\bar{\Upsilon}(-1/z).
\]
\emph{Nota bene}: The definition of $\Upsilon^*$ is just the dimensionally continued Jacobi
transformation [Eq.~\eqref{JTF1}] of $\Upsilon$ with the
factor of $\sqrt{\det\Upsilon}$ omitted. (We leave off this factor, since it
would just cancel against the one present in the
standard Poisson summation formula [cf.\ Eq.~\eqref{nPSF}].)

We now need to show that the generalized theta series have the appropriate properties 
to allow us to copy
the proof of Theorem~\ref{Zd-thm} almost verbatim. Specifically, we need the following
\begin{prop}
If $\Upsilon$ denotes any of the $\Upsilon_d$s or $\Upsilon_d^*$s defined above, we can write
\[
\Upsilon(q) = \sum_{l=0}^\infty N_lq^{A_l},
\]
where
\begin{enumerate}[1.]
\item $A_{l + 1} > A_l$, $A_0 \ge 0$.
\item $\sum_{l=1}^\infty A_l^{-2} < \infty$.
\item There exists $L\in\N$ and $C,n > 0$ such that
$|N_l| \leq C A_l^n$ for all $l \ge L$.
\item The series converges inside the unit disk.
\end{enumerate}
\end{prop}

\begin{proof}
First note that we can apply the same arguments to $\Upsilon_d^*$ as to $\Upsilon_d$, so we
can restrict our attention to the former, without loss of generality. Now, we then have
\[
A_l = (l + \cA)/V, \qquad\qquad \cA := \sum_{m=1}^\cM V\lambda_m s_m/4
\]
where $V$ is the least common
denominator of $s_m$, $t_m$, and $u_m$ (for all $m$). [We have the additive
constant $\cA$ due to the overall factor of $q^{1/4}$ in $\vartheta_2(q)$.]
Thus the first required and second required properties (positivity and monotonicity of the $A_l$
and convergence of the series whose terms are $A_l^{-2}$) are 
obviously true.

For the third property (polynomial boundedness of the $N_l$), we use the same
Cauchy's integral formula argument used in Sec.~\ref{Nl_bounds}. (The analyticity established
below shows that Cauchy's theorem is still applicable here.) Here $N_l$
is given by the $l$th term in the Maclaurin expansion of
$\Upsilon_d(q^V)/q^\cA$, so we have
\<\nonumber
\begin{split}
|N_l| &= \left|\frac{1}{2\pi i}\int_{\mathsf{C}_R} \prod_{m=1}^\cM\frac{\vartheta_2^{\lambda_m}(z^{V{s_m}})\vartheta_3^{\rho_m}(z^{Vt_m})\vartheta_4^{\sigma_m}(z^{Vu_m})}{z^{V\lambda_m s_m/4}z^{l+1}}dz\right|\\
&\le 2^d\prod_{m=1}^\cM\frac{1}{R^l(1 - R^{Vs_m})^{\lambda_m}(1 - R^{Vt_m})^{\rho_m}(1 - R^{Vu_m})^{\sigma_m} }
\le \frac{2^d}{R^l(1-R)^d},
\end{split}
\?
where $\mathsf{C}_R$ is the same contour used previously.
We have used the geometric series to obtain
the bound $|\tau(q)| \le 2/(1-|q|)$, where $\tau(q)$ is any of
$\vartheta_2(q)/q^{1/4}$, $\vartheta_3(q)$, or $\vartheta_4(q)$. 
Additionally, we have used
the fact that $\kappa \ge 1$, where $\kappa$ is any of $Vs_m$, $Vt_m$, or $Vu_m$, so
$|1 - R^\kappa| \ge 1- R$, since $R \in (0,1)$. We also
recalled that $\lambda_m,\rho_m,\sigma_m\geq 0$ and $\sum_{m=1}^\cM(\lambda_m+\rho_m+\sigma_m) = d$. Since there is an $R\in(0,1)$ such that
$2^d/[R^l (1-R)^d] \le C_d l^d$ (for $l \geq 1$),
as was shown in Sec.~\ref{Nl_bounds}, we are done.

The fourth property (convergence of the $q$-series in the unit disk) follows
from the analyticity and lack of zeros of the theta functions inside the unit
disk, as in the argument
given below Eq.~\eqref{Taylor}. [Note that here we consider $\vartheta_2(q)/q^{1/4}$, not $\vartheta_2(q)$ itself.] Specifically, $\Upsilon(q^V)/q^\cA$ is an
analytic function of $q$ inside the unit disk; the lack of zeros
can be seen from the infinite product representations of $\vartheta_2$ and
$\vartheta_4$ given, e.g., in Eqs.~(34) and (36) in Chap.~4 of Conway and
Sloane~\cite{CS}.
\end{proof}

The dimensionally continued summation formula for generalized theta series thus takes the form of the following
\begin{thm}\label{General_theorem}
Let $\Upsilon$ be a generalized theta series (as defined above), with power series coefficients and powers $N_l$ and
$A_l$, i.e.,
\[
\Upsilon(q) = \sum_{l = 0}^\infty N_lq^{A_l}.
\]
Let $N_l^*$ and $A_l^*$ be the corresponding quantities for $\Upsilon$'s dual, $\Upsilon^*$
[computed using Eq.~\eqref{Upsilon_d_star}]. Then, for any $f\in\sS^\sfE$, we have the
summation formula
\[
\sum_{l=0}^\infty N_lf(\sqrt{A_l}) = \sum_{l=0}^\infty N^*_l\hat{f}(\sqrt{A_l^*}),
\]
where we compute $\hat{f}$ using Eq.~\eqref{f-hat} (with the dimension parameter $d$ associated
with $\Upsilon$).
\end{thm}

\begin{rem}
In general, all one requires of the $\Upsilon$s used in this summation formula is that they and their
Jacobi transformations have sufficiently well-behaved power series. (This is satisfied by all functions
of \emph{weight}, in the terminology of Pache~\cite{Pache}---see Pache's Definition~11.) However, it is unclear whether any 
such functions exist besides the generalized theta series we have defined above, so we have not given the theorem in a more general form.
\end{rem}

\begin{proof}
The proof is almost the same as that for Theorem~\ref{Zd-thm} (replacing
$\Theta$ by $\Upsilon$, and noting that we can no longer appeal to the
standard Poisson summation formula for $d = 1$, so we simply exclude that
case). Most of the work has been done in the proof of the Proposition; the only new part is checking that
$\sum_{l=0}^\infty|N_lh(\sqrt{A_l})| \to 0$ as
$\epsilon \to 0$ if $h$ is $\epsilon$-close to $0$ in the Schwartz space
topology [and similarly for $\sum_{l=0}^\infty|N^*_lh(\sqrt{A^*_l})|$]. To do this, we simply note
that we have $|N_lh(\sqrt{A_l})| \le CA_l^n|h(\sqrt{A_l})|$, by polynomial boundedness of the
$N_l$,
and that $x^{2n + 4}|h(x)| \le \epsilon$ $\forall$ $x\in\R$ $\Rightarrow$ $A_l^n|h(\sqrt{A_l})| \le \epsilon/A_l^2$ [cf.\ the discussion at the end of Sec.~\ref{Nl_bounds}], from which the desired result follows immediately.
(The same argument holds for the starred quantities, since they have the same properties as the
unstarred quantities.)
\end{proof}

\begin{rem}
This theorem can likely be interpreted as a trace formula for the dimensionally
continued, spherically symmetric Laplacian [Eq.~\eqref{d-Laplacian}], since the kernel of the
dimensionally continued Fourier transform is an eigenfunction of this operator
(see Lemma~\ref{eig-lem}).
See, e.g., Sec.~1.3 (particularly Theorem~1.3) of Uribe~\cite{Uribe} for a
presentation of the standard Poisson summation formula for an integer
dimension lattice as a trace formula for the Laplacian.
\end{rem}

\begin{rem}
This result shows that one can apply this extended Poisson summation formula to lattice-like
objects whose theta series have coefficients of both signs, so they do not exist as a lattice,
even though $d\in\N$: For a trivial example, consider $d = 2$ and $\Phi(q) = \vartheta_4^2(q) = 1 - 4q + 4q^2 + \cdots$. Of course, this is in some sense
too trivial, since one can write $\vartheta_4^2 = 2\Theta_{D^2} -
\Theta_{\Z^2}$, and then apply the standard Poisson summation formula to each
of those lattices to establish the result in this case (cf.\ the discussion
in C\'{o}rdoba~\cite{Cordoba_LMP}).
However, in more
complicated higher-dimensional cases, it will likely not be clear how to
construct the lattice(s) associated with the theta series (if they indeed
exist). Indeed, Jenkins and Rouse~\cite{JeRo} have very recently shown that for weights higher
than $81\,632$, all the modular forms of a certain type have coefficients of both signs. (This general property was first
shown
in less specific, much earlier work by Mallows, Odlyzko, and Sloane~\cite{MOS}.)
\end{rem}

\section{Outlook}

While Theorem~\ref{General_theorem} encompasses quite a large family of summation formulae, 
there still remains wide latitude for further generalizations (even excluding the various possibilities
for weakening certain of the hypotheses mentioned after Theorem~\ref{Zd-thm}). The most
sweeping generalization would likely be to replace the dimensionally continued Fourier
transform with some more general family of integral transforms, with the possibility of a
subsequent enlargement of the transformation properties required of the generalized theta series.
Here one could follow the work of Kubota	~\cite{Kubota} and Unterberger~\cite{Unterberger} in
integer dimensions. But even if one retains the dimensionally continued Fourier transform, one
can still likely obtain summation formulae from more general classes of generalized theta series
than we have considered. In particular, it would be interesting to obtain a dimensionally continued
version of the quasicrystal summation formula given as Theorem~2.9 in Lagarias~\cite{Lagarias}.
Here the calculations of the central shelling for certain quasicrystals in, e.g.,~\cite{BaGr, MP, MW, MoSa} could be relevant.
Additionally, since the coefficients of standard theta series give the representation numbers for lattices, it is
possible that our results could be applicable to generalizations of representation number
problems: See, e.g.,~\cite{Hanke} for a review of standard results on representation numbers.

\section*{Acknowledgements}

It is our pleasure to thank George Andrews, Jingzhi Tie, and
Shuzhou Wang for encouragement and useful comments, John Roe for
mentioning an alternative proof of the density result,
Artur T\v{s}obanjan for
helpful remarks about the exposition, and the anonymous referee for various
perceptive comments and suggestions. We also thank Ben Owen for comments
on the manuscript and for suggesting the physics problem that sparked these 
investigations.
This work was supported 
by NSF grants PHY-0555628 and PHY-0855589, the Eberly research funds of Penn State, and the DFG SFB/Transregio 7.

\bibliographystyle{amsplain}
\bibliography{PSF}

\end{document}